\documentclass[11pt,leqno,letterpaper]{amsart}

\usepackage{amsthm}
\usepackage{mathrsfs}
\usepackage{amsmath}
\usepackage{amsmath,amsthm,amssymb,amscd}
\usepackage{latexsym}
\usepackage{enumerate}
\usepackage[utf8]{inputenc}
\usepackage[english]{babel}

\usepackage{graphicx}
\usepackage{color}
\usepackage{bm} 
\usepackage{amsfonts} 
\usepackage{multirow}
\usepackage{tikz}

\newtheorem{theorem}{Theorem}[section]
\newtheorem{proposition}{Proposition}[section]
\newtheorem{lemma}{Lemma}[section]
\newtheorem{corollary}{Corollary}[section]

\newtheorem{remark}{Remark}[section]

\newcommand{\tr}{^{\sf T}}
\newcommand{\m}[1]{{\bf{#1}}}
\newcommand{\g}[1]{\bm #1}
\newcommand{\C}[1]{{\mathcal {#1}}}

\newcommand{\Supp}[1]{{\rm supp}\; (#1)}

\usepackage[colorlinks=true,breaklinks=true,bookmarks=true,urlcolor=blue,
     citecolor=blue,linkcolor=blue,bookmarksopen=false,draft=false]{hyperref}

\begin{document}

\title[A General Regularized Continuous Formulation for the MCP]{A General Regularized Continuous Formulation for the Maximum Clique Problem}

\author[J. T. Hungerford, F. Rinaldi]{James T. Hungerford$^1$ \MakeLowercase{and} Francesco Rinaldi$^2$}

\address{RaceTrac Store Support Center\\ 200 Galleria Pkwy SE\\ Atlanta, GA 30339}

\email{jamesthungerford@gmail.com}

\address{Dipartimento di Matematica\\Universit\`{a} di Padova\\Via Trieste 63\\35121 Padova - Italy}

\email{rinaldi@math.unipd.it}



\begin{abstract}
In this paper, we develop a general regularization-based continuous optimization framework 
for the maximum clique problem. In particular, we consider a broad class of regularization terms
that can be included in the classic Motzkin-Strauss formulation
and we develop conditions that guarantee the equivalence between the continuous regularized problem 
and the original one in both a global and a local sense. 
We further analyze, from a computational point of view, two different regularizers that satisfy the general conditions.
\end{abstract}

\maketitle

\section{Introduction.}
\label{introduction}
Let $G = (\C{V},\C{E})$ be a simple undirected graph on vertex set
$\C{V} = \{1,2,\ldots,n\}$ and edge set $\C{E}\subseteq\C{V}\times\C{V}$.
Since $G$ is simple and undirected
$(j,i)\in\C{E}$ whenever $(i,j)\in\C{E}$,
and $(i,i)\notin\C{E}$ for
any $i\in\C{V}$.
A \emph{clique} in $G$ is
a subset $C \subseteq \C{V}$ such that $(i,j)\in\C{E}$ for every $i,j\in C$
with $i\neq j$. In this paper, we consider the
classical \emph{Maximum Clique Problem} (MCP):
find a clique $C\subseteq \C{V}$ such that $|C|$ is maximum.

The Maximum Clique Problem has a wide range of applications
(see \cite{Bomze1999, Wu2015}
and references therein) in areas such as
social network analysis, telecommunication networks,
biochemistry, and scheduling.
The cardinality of a maximum clique in $G$ is denoted $\omega (G)$.
A clique $C$ is said to be \emph{maximal} if it is not contained in any
strictly larger clique; that is, if there does not exist a clique
$D$ such that
$C\subset D$. $C$ is said to be \emph{strictly maximal} if there do not
exist vertices $i\in C$ and $j\notin C$ such that $C\cup\{j\}\backslash\{i\}$
is a clique.

The MCP is NP-hard \cite{Karp72}. However, due in part
to its wide applicability,
a large variety of both heuristic and
exact approaches have been investigated
(see \cite{Bomze1999}
for a thorough overview of formulations and algorithms going up to 1999; a 
more recent survey of algorithms is given in \cite{Wu2015}). A significant
number of the solution methods proposed
(for example,
\cite{bomze1997evolution,Bomze2002,Gibbons97,Kuznetsova2001,Motzkin65, Pelillo1996, Pelillo1996b})
are based on solving the following
well-known
continuous quadratic programming formulation of the MCP, due to Motzkin
and Straus
\cite{Motzkin65}:
\begin{eqnarray}\label{MSQP}
&\max\quad\m{x}\tr\m{A}\m{x}\\
&\mbox{subject to}\quad\m{x}\in\Delta\;, \nonumber
\end{eqnarray}
where
$\Delta$ is the $n$-dimensional simplex defined by
\[
\Delta := \{\m{x}\in\mathbb{R}^n : \m{0}\le\m{x}\le\m{1} \mbox{  and  }
\m{1}\tr\m{x} = 1\}
\]
and $\m{A} = (a_{ij})_{i,j\in\C{V}}$
denotes the adjacency matrix for $G$ defined by
\[
a_{ij} =
\left\{
\begin{array}{rl}
1,& (i,j)\in\C{E}\\
0,& (i,j)\notin\C{E}
\end{array}
\right.\quad\forall\; i,j\in\C{V}\;.
\]
For any non-empty
clique $C$ we let $\m{x} (C)\in\Delta$ denote
the corresponding characteristic vector (defined by
$x (C)_i = \frac{1}{|C|}$ whenever $i\in C$ and $x (C)_i = 0$ otherwise).
The equivalence between MCP and (\ref{MSQP}) is given by
the following theorem:
\smallskip

\begin{theorem}[Theorem 1 \cite{Motzkin65}]\label{MStheorem} 
The optimal objective value of {\rm (\ref{MSQP})} is
\[
1 - \frac{1}{\omega(G)}
\]
and
$\m{x} (C)$ is a global maximizer of {\rm (\ref{MSQP})}
for any maximum clique $C$.
\end{theorem}
\smallskip

Solution approaches to MCP based on solving (\ref{MSQP}) include
nonlinear programming
methods \cite{Gibbons93} and methods based on
discrete time
replicator dynamics
\cite{bomze1997evolution, Bomze1999, Bomze2002, Pelillo1996}.
Since (\ref{MSQP}) is NP-hard (by reduction to MCP),
the computing time required to obtain a global maximizer
can grow exponentially with the size of the graph; hence,
finding a global maximizer may be impractical in many settings.
On the other hand, iterative optimization methods
will typically converge to a point satisfying the first-order
optimality (Karush-Kuhn-Tucker) conditions.
In general, verifiying whether a first-order point of a quadratic
program is even \emph{locally} optimal is an NP-hard problem
\cite{murty1987,pardalos88}.
However, it was shown
in \cite{Gibbons97}
that local optimality of a first-order point
(in fact, any feasible point)
in (\ref{MSQP})
can be ascertained in polynomial time.

In \cite[Proposition 3]{Pelillo1996b},
a characteristic vector for a clique was shown to satisfy
the standard first-order optimality condition for (\ref{MSQP}) if
and only if the associated clique is maximal.
In \cite[Theorem 2]{Gibbons97}, the authors gave a characterization of the
local optima of (\ref{MSQP}) and demonstrated a one-one correspondence
between strict local maximizers 
and strictly maximal cliques.
These results suggest the
possibility of applying iterative optimization methods to (\ref{MSQP})
in order to \emph{approximately} solve MCP (ie. to find large \emph{maximal}
cliques).
However, one known \cite{bomze1997evolution, Pelillo1996, Pelillo1996b}
drawback of this approach in practice
is the presence of ``infeasible'' or ``spurious''
local maximizers of (\ref{MSQP})
which are not 
characteristic vectors for cliques and from which a clique can
not be recovered through any simple transformation.
Such points are an undesirable property of the program,
since they can cause continuous based heuristics to fail by terminating
without producing a clique.
In \cite{bomze1997evolution},
the author addresses this issue by introducing the following
regularized formulation (with $\alpha = \frac{1}{2}$):
\begin{eqnarray}\label{BQP}
&\max\quad \m{x}\tr\m{A}\m{x} + \alpha ||\m{x}||_2^2\\
&\mbox{subject to} \quad \m{x}\in\Delta\;. \nonumber
\end{eqnarray}
In contrast to (\ref{MSQP}),
the local maximizers of (\ref{BQP}) have been shown to be
in one-one correspondence
with the
maximal cliques in $G$ (see \cite[Theorem 9]{bomze1997evolution}), and
a replicator dynamics approach to solving (\ref{BQP}) was shown to reduce
the total number of algorithm failures by 30\%,
compared with 
a similar approach to solving (\ref{MSQP}).
In \cite{Bomze2002}, the authors
enhanced the algorithm of \cite{bomze1997evolution}, adding an annealing
heuristic to obtain even stronger results. In addition, it was demonstrated
that the correspondence between the local/global optima of (\ref{BQP})
and MCP is maintained for any $\alpha \in (0,1)$.
A similar formulation and approach \cite{Bomze2000}
has also
been applied successfully to a weighted version of MCP.
In practice, the numerical performance of 
an
iterative optimization method (in terms of speed and/or solution quality)
may depend on the particular regularization
term employed. 
In this paper we consider a broad class of regularization terms, and
develop conditions
under which
the regularized program is equivalent to MCP in both a global and a local
sense. We establish the equivalence in a step by step manner that reveals
some of the underlying structural properties
of (\ref{MSQP}).
We provide two different examples of regularization
terms satisfying the general conditions, 
and also give
some preliminary computational results evaluating their
effectiveness in terms of both speed and solution quality.
Over the course of our analysis (see Section \ref{sectGenReg}),
we also correct an (apparently as yet
unidentified) erroneous result in the literature
linking maximal cliques to
local maximizers in (\ref{MSQP}) by constructing
an example of a maximal clique whose characteristic vector
is not a local maximizer of (\ref{MSQP}).

The paper is organized as follows. In Section 2,
we develop a general regularized formulation of MCP and provide
conditions under which the global/local maximizers of the regularized
program are in one-one correspondence with the maximum/maximal cliques in $G$.
In Section 3, we report on some preliminary computational results
comparing the performance of two new regularization terms with the one proposed by Bomze in \cite{bomze1997evolution}.
We conclude in Section 4.

{\bf Notation.}
\m{0} and \m{1} denote column vectors whose entries are all 0 and all 1 respectively and
$\m{I}$ denotes the
identity matrix,
where the 
dimensions should be clear from the context.
$\nabla f (\m{x})$ denotes the gradient of $f$, a row vector,
and $\nabla^2 f(\m{x})$ denotes the Hessian of $f$.
For a set $\C{Z}$, $|\C{Z}|$ is the number of elements in $\C{Z}$.
$\m{e}_i\in\mathbb{R}^n$
denotes the $i$-th column of the $n\times n$ identity matrix.
If $\{s_i\}_{i = 1}^n
\subset \mathbb{R}$
is a finite sequence of length $n$, then $Diag (\{s_i\}_{i = 1}^n)$ is the
$n\times n$ diagonal matrix whose $(i,i)$th entry is $s_i$.
If $\m{x}\in\mathbb{R}^n$, then
$\Supp{\m{x}}$ denotes the \emph{support} of $\m{x}$,
defined by
$\Supp{\m{x}} = \{ i : x_i \ne 0 \}$.
Given vectors $\m{x},\m{y}\in\mathbb{R}^n$, $[\m{x},\m{y}] := \{t\m{x} + (1-t)\m{y}\; :\; t\in [0,1]\}$.
For a given positive integer $n$,
we denote the set $\{1,2,\ldots,n\}$ by $[n]$.
$\C{S}_n$ is the set of permutations of $[n]$.
If $\m{B}\in\mathbb{R}^{n\times n}$ is a symmetric matrix, we write
$\m{B}\succeq \m{0}$ if $\m{B}$ is positive semidefinite, $\m{B}\succ \m{0}$
if $\m{B}$ is positive definite, and $\m{B} \preceq \m{0}$
(resp. $\m{B} \prec \m{0}$) when $-\m{B} \succeq \m{0}$
(resp. $-\m{B}\succ \m{0}$).
If $\m{X}\subseteq\mathbb{R}^n$ and $f:\m{X} \rightarrow \mathbb{R}$,
then a point
$\m{x}\in \m{X}$ is a \emph{local} (resp. \emph{strict local})
\emph{maximizer} 
of the problem $\max\;\{f (\m{x}) \;:\;\m{x}\in \m{X}\}$ if
there exists some $\epsilon > 0$ such that
$f (\m{x}) \ge f (\tilde{\m{x}})$ (resp.
$f (\m{x}) > f (\tilde{\m{x}})$)
for every
$\tilde{\m{x}}\in\m{X}$
with $0 < ||\tilde{\m{x}} - \m{x}||_2 < \epsilon$.
$\m{x}$ is an \emph{isolated} local maximizer if there exists some
$\epsilon > 0$ such that $\tilde{\m{x}}$ is not a local maximizer for any
$\tilde{\m{x}}\in\m{X}$ with $0 < ||\tilde{\m{x}} - \m{x}||_2 < \epsilon$.
$conv(\m{X})$ denotes the convex
hull of $\m{X}$. $span_+(\m{X}) := \{\sum_{i = 1}^k \alpha_i \m{x}^i\;:\; k\in\mathbb{N}_{+},\;\alpha_i\ge 0,\;\m{x}^i\in \m{X}\;\forall\;i\in [k]\}$. $\mathbb{N}_{+}$ denotes the
set of positive natural numbers.

\section{General regularized formulation.}
\label{sectGenReg}

Consider the following problem:
\begin{eqnarray}\label{PMS}
&\max\quad f (\m{x}) := \m{x}\tr\m{Ax} + \Phi (\m{x})\\
&\mbox{subject to} \quad \m{x}\in\Delta\;, \nonumber
\end{eqnarray}
where $\Phi : X \rightarrow \mathbb{R}$ is a twice continuously differentiable
function defined on some open set
$X \supset \Delta$.
Throughout this section, we will also assume
that
$\Phi$ satisfies the following conditions for every $\m{x}\in\Delta$:
\smallskip

\begin{itemize}
\item[(C1)] $\nabla^2 \Phi (\m{x}) \succeq \m{0}$
\item[(C2)] $||\nabla^2\Phi (\m{x})||_2 < 2$
\item[(C3)] $\Phi$ is constant on the set
\[
\C{P} (\m{x}) :=
\{\tilde{\m{x}}\in\Delta\;:\;\exists\;\sigma\in\C{S}_n\mbox{ such that } \tilde{x}_i = x_{\sigma (i)}\;\forall\;i\in [n]\}\;,
\]
\end{itemize}
\smallskip
\noindent
where $\C{S}_n$ is the set of permutations of $[n]$.
Note that (C1) is equivalent to requiring that $\Phi$ is
convex at $\m{x}$ and
(C2) is equivalent to
$\nabla^2 \Phi (\m{x}) - 2\m{I} \prec \m{0}$
(a fact that will be used later).
Also,
since $\Phi \equiv 0$
satisfies (C1) -- (C3) trivially,
the results
of this section will hold in particular
for the original unpenalized formulation (\ref{MSQP}) when no
additional assumptions are made on $\Phi$.

We will establish the global equivalence between (\ref{PMS})
and MCP through a series of intermediate results.
For any clique $C$ define the set
\[
\Delta (C) := \{\m{x}\in\Delta\; : \; \Supp{\m{x}} \subseteq C\}\;,
\]
and let
\begin{equation}\label{Delta0def}
\displaystyle{\Delta^0 \;\; := \;\; \bigcup_{C \mbox{ clique}} \Delta (C)\;\; = \;\;
\{\m{x}\in\Delta\;:\;\Supp{\m{x}}\mbox{ is a clique}\}}\;.
\end{equation}

\begin{lemma}\label{lemma1}
Let $\m{x}\in\Delta^0$. Then,
\begin{itemize}
\item [{\rm 1.}] For any
$\tilde{\m{x}}\in\C{P} (\m{x})\cap\Delta^0$ we have
\[
f (\tilde{\m{x}}) = f (\m{x})\;.
\]
\item [{\rm 2.}] For any $\m{0}\neq\m{d}\in\mathbb{R}^n$ such that
$\m{x} + t\m{d}\in\Delta^0$ for all sufficiently small $t > 0$ we have
\[
\m{d}\tr\nabla^2 f (\m{x}) \m{d} < 0\;.
\]
\end{itemize}
\end{lemma}
\smallskip

\proof
First, we note that
for any $\m{z}\in\mathbb{R}^n$ such that $\Supp{\m{z}}$ is
a clique, we have
\begin{eqnarray}
\m{z}\tr\m{A}\m{z} & = &
\sum_{i = 1}^n \sum_{j = 1}^n z_i a_{ij} z_j \;=\;\sum_{i\in\Supp{\m{z}}}\sum_{j\in\Supp{\m{z}}\backslash\{i\}} z_iz_j\nonumber\\
& = &
\sum_{i\in\Supp{\m{z}}} \left(\sum_{j\in\Supp{\m{z}}} z_i z_j - z_i^2\right)\nonumber\\
\\
& = &
(\m{1}\tr\m{z})^2 - \m{z}\tr\m{z}\;.\label{eq20}
\end{eqnarray}
Next let $\m{x}\in\Delta^0$. We prove the two parts separately.

{\bf Part 1.}
For any $\tilde{\m{x}}\in\C{P} (\m{x})\cap\Delta^0$ we have
\begin{eqnarray}
f (\tilde{\m{x}}) & = &
\tilde{\m{x}}\tr\m{A}\tilde{\m{x}} + \Phi (\tilde{\m{x}})\nonumber\\
& = &
(\m{1}\tr\tilde{\m{x}})^2 - \tilde{\m{x}}\tr\tilde{\m{x}} + \Phi (\tilde{\m{x}})\label{eq40}\\
& = &
(\m{1}\tr{\m{x}})^2 - {\m{x}}\tr{\m{x}} + \Phi ({\m{x}})\label{eq41}\\
& = & {\m{x}}\tr\m{A}{\m{x}} + \Phi ({\m{x}})\label{eq42}\\
& = & f (\m{x})\;,\nonumber
\end{eqnarray}
where (\ref{eq40}) and (\ref{eq42}) are due to (\ref{eq20}), since
$\tilde{\m{x}},\m{x}\in\Delta^0$, and (\ref{eq41}) is due to {\rm(C3)} and the
assumption that $\tilde{\m{x}}\in\C{P} (\m{x})$.

{\bf Part 2.}
Let $\m{0}\neq\m{d}\in\mathbb{R}^n$ be any vector such that $\m{x} + t \m{d} \in \Delta^0$ for all sufficiently small $t > 0$. Then when $t$ is sufficiently
small $\Supp{\m{d}}\subseteq\Supp{\m{x} + t\m{d}}$, which implies
that $\Supp{\m{d}}$ is a clique; moreover, since
$\m{1}\tr\m{x} = \m{1}\tr(\m{x} + t \m{d}) = 1$, we have $\m{1}\tr\m{d} = 0$.
So by (\ref{eq20}) we have
\begin{eqnarray}
\m{d}\tr\nabla^2 f (\m{x}) \m{d} &\; =\;&2\m{d}\tr\m{A}\m{d} + \m{d}\nabla^2\Phi (\m{x})\m{d}\\
&\;=\; & 2 [(\m{1}\tr\m{d})^2 - \m{d}\tr\m{d}] +
\m{d}\tr\nabla^2 \Phi (\m{x})\m{d}\nonumber\\
&\; =\; &
-2\m{d}\tr\m{d} + \m{d}\tr\nabla^2\Phi (\m{x})\m{d}\nonumber\\
&\; =\; &
-\m{d}\tr\left[2\m{I} - \nabla^2 \Phi (\m{x})\right]\m{d}\nonumber\\
&\; < \;& 0\;,\label{inequality1}
\end{eqnarray}
where (\ref{inequality1}) follows from {\rm(C2)}.
\endproof
\smallskip

Now consider the following problem:
\begin{eqnarray}\label{PDeltaC}
& \max\quad f (\m{x})\\
& \mbox{subject to}\quad \m{x}\in\Delta (C)\;.\nonumber
\end{eqnarray}

\begin{proposition}\label{lemmaDeltaC}
The unique local (hence global) maximizer of
{\rm (\ref{PDeltaC})} is $\m{x} (C)$.
\end{proposition}
\smallskip

\proof
Suppose by way of contradiction that there exist distinct local
maximizers $\m{x}^1 \neq \m{x}^2$ of (\ref{PDeltaC}). Then by Taylor's Theorem
(see for example Proposition A.23 of \cite{be99})
$\exists\;\g{\xi}\in[\m{x}^1,\m{x}^2]\subseteq\Delta$ such that
\begin{equation}\label{taylor1}
f (\m{x}^2) = f (\m{x}^1)  + \nabla f (\m{x}^1)\m{d} + \frac{1}{2}\m{d}\tr
\nabla^2 f (\g{\xi})\m{d}\;,
\end{equation}
where $\m{d} = \m{x}^2 - \m{x}^1$.
Since $\m{x}^1,\m{x}^2\in\Delta$, we have that $\m{x}^1 + t\m{d}\in[\m{x}^1,\m{x}^2]\subseteq\Delta (C)$ for all sufficiently small $t > 0$. So by the standard
first-order necessary local optimality condition
(see for example, Section 1 of
\cite{HagerHungerford13}), we have
\begin{equation}\label{firstorderdescent}
\nabla f (\m{x}^1)\m{d} \le 0\;.
\end{equation}
Moreover, Part 2 of Lemma \ref{lemma1} implies
\begin{equation}\label{suffcond}
\m{d}\tr\nabla^2 f (\g{\xi})\m{d} < 0\;.
\end{equation}
Combining (\ref{taylor1}),
(\ref{firstorderdescent}), and (\ref{suffcond}), we obtain
$f (\m{x}^2) < f (\m{x}^1)$. But then
interchanging $\m{x}^2$ and $\m{x}^1$ the
same argument can be used to show that $f (\m{x}^1) < f (\m{x}^2)$,
a contradiction. Therefore, there is a unique local (hence global) maximizer
of (\ref{PDeltaC}), say $\m{x}^*$.

Next, we claim that $\C{P} (\m{x}^*)\cap\Delta(C) = \{\m{x}^*\}$.
Indeed, suppose by way of contradiction that $\exists\;\tilde{\m{x}}\in\C{P} (\m{x}^*)\cap\Delta (C)$ such that $\tilde{\m{x}} \neq \m{x}^*$. Since
$\Supp{\m{x}^*}\subseteq C$ and $C$ is a clique, Part 1 of
Lemma \ref{lemma1} implies
$f (\tilde{\m{x}}) = f (\m{x}^*)$. But then $\tilde{\m{x}}$ must be a global
maximizer of (\ref{PDeltaC}), contradicting the uniqueness of $\m{x}^*$. Hence,
we must have $\C{P} (\m{x}^*)\cap\Delta (C) = \{\m{x}^*\}$.
Thus, $x^*_i = x^*_j$ for any $i,j\in C$. Since
$\m{x}^*\in\Delta (C)$, this implies $\m{x}^* = \m{x} (C)$.
\endproof
\smallskip

\begin{remark}
By Part 2 of Lemma \ref{lemma1}, 
{\rm(\ref{PDeltaC})} is a
strictly concave (and smooth) maximization problem. Thus, the uniqueness
of the maximizer of {\rm(\ref{PDeltaC})}
may be seen as following from standard results in the theory
of convex optimization
(for instance, Proposition B.4 in \cite{be99}).
\end{remark}
\smallskip

Next, consider the problem
\begin{eqnarray}\label{Delta0P}
& \max\quad f (\m{x})\\
& \mbox{subject to}\quad \m{x}\in\Delta^0\;.\nonumber
\end{eqnarray}

\begin{proposition}\label{lemmaDelta0P}
A point $\m{x}\in\Delta^0$ is a local maximizer of {\rm (\ref{Delta0P})}
if and only if $\m{x} = \m{x} (C)$ for some maximal clique $C$. Moreover,
every local maximizer of {\rm (\ref{Delta0P})} is strict.
\end{proposition}
\smallskip

\proof
First,
observe that for any local maximizer $\m{x}$ of 
(\ref{Delta0P}), by (\ref{Delta0def})
there exists some maximal clique $C$
such that $\m{x}\in\Delta (C)$; and since $\m{x}$ is a local maximizer
of (\ref{Delta0P}), it is also a local maximizer of (\ref{PDeltaC}),
which implies
that
$\m{x} = \m{x} (C)$, by Proposition \ref{lemmaDeltaC}.
Thus,
the proof will be complete when we show that every characteristic vector for
a maximal clique is a strict local maximizer in (\ref{Delta0P}).
To this end,
let $C$ be a maximal clique, and suppose by way of contradiction that
$\m{x} (C)$ is not a strict local maximizer of (\ref{Delta0P}). Then,
for every $k\in\mathbb{N}_{+}$ there exists some $\m{x}^k\in\Delta^0$
with $0 < ||\m{x}^k - \m{x} (C)||_2 < 1/k$
such that
$f (\m{x}^k) \ge f (\m{x} (C))$.
Since there are only finitely many sets in the unions in (\ref{Delta0def}),
there must exist some clique $C'$ and some subsequence
$(\m{x}^{k_l})_{l = 1}^{\infty} \subseteq (\m{x}^k)_{k = 1}^{\infty}$ such
that $\m{x}^{k_l}\in\Delta (C')$ for each $l \ge 1$, with
$\m{x}^{k_l}\rightarrow\m{x} (C)$. Hence,
$\m{x} (C)\in\overline{\Delta(C')} = \Delta (C')$, which implies
$C = \Supp{\m{x} (C)} \subseteq C'$. Since $C$ is maximal, we must
have that $C = C'$, and thus $\m{x}^{k_l}\in\Delta (C') = \Delta (C)$
for each $l \ge 1$. Thus, $\m{x} (C)$ is not a strict local maximizer
of (\ref{PDeltaC}),
contradicting Proposition \ref{lemmaDeltaC}.
This completes the proof.
\endproof
\smallskip

\begin{proposition}\label{propcliquesize}
If $C^1$ and $C^2$ are cliques, then
\[
|C^1| < |C^2| \;\Leftrightarrow\;
f ( \m{x} (C^1) ) < f ( \m{x} (C^2) )\;.
\]
\end{proposition}
\smallskip

\proof
Let $C^1$ and $C^2$ be cliques.
First, suppose that 
$|C^1| < |C^2|$.
Let $C$ be any clique such that $C \subset C^2$ and $|C| = |C^1|$. Then
$\m{x} (C^1)\in\C{P} (\m{x} (C))$. So, Part 1 of Lemma \ref{lemma1} implies
$f (\m{x} (C^1)) = f (\m{x} (C))$. Moreover, by Proposition \ref{lemmaDeltaC}
$f (\m{x} (C)) < f (\m{x} (C^2))$, since $\m{x} (C)\in\Delta (C^2)$.
Hence, $f (\m{x} (C^1)) < f (\m{x} (C^2))$.
Conversely, suppose that
$f (\m{x} (C^1)) < f (\m{x} (C^2))$.
Then, by the proof of the forward direction we must have $|C^1| \le |C^2|$.
Moreover, if
$|C^1| = |C^2|$, then $\m{x} (C^1) \in\C{P} (\m{x} (C^2))$ and
Part 1 of Lemma \ref{lemma1} implies $f (\m{x} (C^1)) = f (\m{x} (C^2))$, a
contradiction. Hence, we must have $|C^1| < |C^2|$.
\endproof
\smallskip

\begin{corollary}\label{corollaryDelta0P}
A point
$\m{x}\in\Delta^0$ is a global maximizer of {\rm (\ref{Delta0P})}
if and only if
$\m{x} = \m{x} (C)$ for some maximum clique $C$.
\end{corollary}
\smallskip

\proof
Let $\m{x}\in \Delta^0$. Then $\m{x}$ is a global maximizer of
(\ref{Delta0P}) if and only if
$\m{x}$ is a local maximizer and
$f (\m{x}) \ge f (\bar{\m{x}})$ for every
local maximizer $\bar{\m{x}}\neq\m{x}$, which by
Proposition \ref{lemmaDelta0P}
holds if and only if
$\m{x} = \m{x} (C)$ for some maximal clique $C$ and
$f (\m{x} (C)) \ge f (\m{x} (\bar{C}))$
for every maximal clique $\bar{C}\neq C$. The corollary then follows from
Proposition \ref{propcliquesize}.
\endproof
\smallskip

\begin{proposition}\label{propGlobal}
For every clique $C$, 
$\m{x} (C)$ is a global maximizer of {\rm (\ref{PMS})}
if and only if $C$ is a maximum clique.
\end{proposition}
\smallskip

\proof
We will show that
there exists a global maximizer of (\ref{PMS}) which lies in $\Delta^0$.
The proof will then
follow from Corollary \ref{corollaryDelta0P}.
To this end,
let $\m{x}$ be any global maximizer of (\ref{PMS}). If
$\m{x}\in\Delta^0$, then we are done.
So, suppose instead that
$\m{x}\notin\Delta^0$.
Then $\Supp{\m{x}}$ is not a clique and there exist indices
$i\neq j\in\Supp{\m{x}}$
such that $a_{ij} = 0$. Next, for any $t \in [-x_i, x_j]$, let
$\m{x} (t) := \m{x} + t(\m{e}_i - \m{e}_j)$,
and observe that
by Taylor's Theorem there exists some
$\g{\xi} \in [\m{x}, \m{x}(t)]$ such that
\begin{eqnarray}
f (\m{x} (t)) & = &
f (\m{x}) + t \nabla f (\m{x})(\m{e}_i - \m{e}_j) +
\frac{t^2}{2}(\m{e}_i - \m{e}_j)\tr\nabla^2 f (\g{\xi}) (\m{e}_i - \m{e}_j)\nonumber\\
& = &
f (\m{x}) +
\frac{t^2}{2}
\left[
(2a_{ii} + 2a_{jj} - 4 a_{ij}) +
(\m{e}_i - \m{e}_j)\tr\nabla^2\Phi(\g{\xi})(\m{e}_i - \m{e}_j)
\right]\label{byfirstorder}\\
& = &
f (\m{x}) +
\frac{t^2}{2}
\left[
(\m{e}_i - \m{e}_j)\tr\nabla^2\Phi(\g{\xi})(\m{e}_i - \m{e}_j)
\right]\label{boilsdowntoPhi}\\
& \ge & f (\m{x})\;. \label{gxtgreaterthangx}
\end{eqnarray}
Here,  (\ref{byfirstorder}) follows from the first-order optimality
condition at $\m{x}$, which implies
that $\nabla f (\m{x}) (\m{e}_i - \m{e}_j) = 0$
since $\m{x} (t)$ is feasible for all $t\in [-x_i, x_j]$.
Equality (\ref{boilsdowntoPhi}) follows from the fact that
$a_{ij} = 0 = a_{ii} = a_{jj}$. And
(\ref{gxtgreaterthangx})
follows from (C1).
Thus, setting $t = x_j$, we obtain another global maximizer 
$\m{x} (t)\in\Delta$ such that
$\Supp{\m{x} (t)} = \Supp{\m{x}}\backslash\{j\} \subset \Supp{\m{x}}$.
We may repeat this
process, gradually reducing the size of $\Supp{\m{x}}$ while maintaining
global maximality, until
$\Supp{\m{x}}$ is a clique (possibly of size $1$), at which
point the proof is complete, since then $\m{x}\in\Delta^0$.
\endproof

By Proposition \ref{lemmaDelta0P},
a one-one correspondence exists between the local maximizers of
(\ref{Delta0P}) and the maximal cliques in $G$.
However, as is already well-known
in the case when $\Phi\equiv 0$
(see the discussions 
pertaining to ``infeasible'' or ``spurious'' local optima in
\cite{bomze1997evolution, Pelillo1996, Pelillo1996b}),
if $\Delta^0$ in (\ref{Delta0P})
is relaxed to $\Delta \supseteq \Delta^0$ (in fact, $\Delta = conv (\Delta^0)$,
but we need not prove this here), there may exist local
maximizers of (\ref{PMS}) that are \emph{not} characteristic vectors
for cliques,
which may cause iterative optimization methods for
solving (\ref{PMS}) to fail, terminating without producing a clique.

Conversely, when 
$\Phi \equiv 0$
there may exist characteristic vectors for maximal cliques
which are not local maximizers in (\ref{PMS}).
Indeed, in the
graph
$G$ in Figure \ref{cx},
the sets $C = \{1,2\}$
and $\hat{C} = \{3,4,5\}$ are both maximal cliques;
and since
$\m{x} (C),\; \m{x} (\hat{C}) \in \Delta$, we have that for
all sufficiently small $t > 0$, $\m{x} (C) + t \m{d} \in \Delta$, where
$\m{d} = \m{x} (\hat{C}) - \m{x} (C)$.
Moreover, it is easy to check that
by computation one has that $\nabla f (\m{x}) \m{d} = 0$ and
$\m{d}\tr\nabla^2 f (\m{x})\m{d} = \frac{1}{6} > 0$, where $\m{x} = \m{x} (C)$.
Hence, for all 
sufficiently small $t > 0$
\begin{eqnarray*}
f (\m{x} + t \m{d}) & = &
f (\m{x})
+ t \nabla f (\m{x}) \m{d}
+ \frac{t^2}{2} \m{d}\tr\nabla^2 f (\m{x})\m{d}\\
& = &
f (\m{x}) + \frac{t^2}{2}\m{d}\tr\nabla^2 f (\m{x})\m{d}\\
& > & f (\m{x}) \; .
\end{eqnarray*}
Thus, $\m{x} (C)$
is not a local maximizer of (\ref{PMS}), despite the fact that
$C$ is maximal. We note that the above is a counterexample to
\cite[Corollary 2]{Gibbons97}.

\tikzstyle{every node}=[circle, draw, fill=black!50, inner sep=0pt, minimum width = 4pt]
\tikzset{dashed lines/.style={dashed}}
\begin{figure}
\centering
\begin{tikzpicture}
\draw (0, 2) -- (2, 2);
\draw (0, 2) -- (0, 0);
\draw (0, 2) -- (1, 1);
\draw (0, 0) -- (1, 1);
\draw (1, 1) -- (2, 0);
\draw (0, 0) -- (2, 0);
\draw (2, 0) -- (2, 2);
\draw (0, 2) node{};       
\draw (2, 2) node{};       
\draw (1, 1) node{};  
\draw (0, 0)  node{};      
\draw (2, 0)  node{};      
\draw[dashed lines] (1, 2) ellipse [x radius=1.5cm, y radius=0.6cm];
\draw[dashed lines] (1, 0.28) ellipse [x radius=1.54cm, y radius=.96cm];
\draw (0, 2)  node[above=0.1cm,draw=white,fill=white] {${1}$};
\draw (2, 2)  node[above=0.1cm,draw=white,fill=white] {${2}$};
\draw (1, 1)  node[right=0.1cm,draw=white,fill=white] {${3}$};
\draw (0, 0)  node[below=0.1cm,draw=white,fill=white] {${4}$};
\draw (2, 0)  node[below=0.1cm,draw=white,fill=white] {${5}$};
\draw (2, 2) node[right=0.6cm,draw=white,fill=white] {${C}$};
\draw (2, 0.4) node[right=0.6cm,draw=white,fill=white] {${\hat{C}}$};
\end{tikzpicture}
\caption{\label{cx}An example of a graph $G$ and two maximal cliques $C$ and $\hat{C}$. Here, $\m{x} (C)$ is 
not a local maximizer of {\rm(\ref{MSQP})} even though 
$C$ is maximal.}
\end{figure}
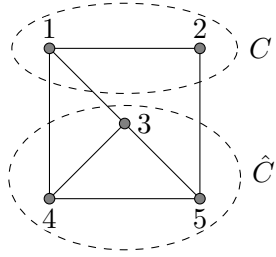

Hence, there is
not necessarily any relationship
(in either direction) between the
local optima of (\ref{MSQP}) and the maximal cliques in $G$.
However,
we will see in the next proposition that
for any strictly convex
$\Phi$ satisfying {\rm(C1)} -- {\rm(C3)},
the
local maximizers of (\ref{PMS}) are in \emph{one-one} correspondence with
the characteristic vectors
for maximal cliques. The proposition is based on three lemmas.
\smallskip

\begin{lemma}\label{lemmaEdgeDesc}
Let $\m{x}\in\Delta$ and let $\C{F}(\m{x})$ denote the set of
first-order feasible directions for {\rm(\ref{PMS})} at $\m{x}$, defined by
\[
\C{F} (\m{x}) = \{\m{d}\in\mathbb{R}^n \; :\; \m{1}\tr\m{d} = 0\mbox{ and }
d_i \ge 0\mbox{ whenever } x_i = 0 \}\;.
\]
Then,
\[
\C{F}(\m{x}) = span_+ (\C{F}(\m{x})\cap\C{D})\;,
\]
where $\C{D} = \bigcup_{\stackrel{i,j = 1}{i\neq j}}^n \{\m{e}_i - \m{e}_j\}$.
\end{lemma}
\smallskip

\proof
The lemma follows immediately from
\cite[Corollary 2.2]{HagerHungerford13} and
the fact that $\C{D}$ is a reflective edge-description
(defined in \cite{HagerHungerford13})
of $\Delta$.
\endproof
\smallskip

The next
lemma is a restatement in the language of the present paper
of the result
\cite[Proposition 3]{Pelillo1996b}
which states that a characteristic vector for a maximal clique
satisfies the first-order optimality conditions of (\ref{MSQP}).
\smallskip

\begin{lemma}\label{firstorderLemma}
If $C$ is a maximal clique, then
\[
\m{x} (C)\tr\m{A} \m{d} \le 0 \quad \forall\;\m{d}\in\C{F}(\m{x} (C))\;.
\]
\end{lemma}

\proof
Let $C$ be a maximal clique. 
We claim that for any
$\m{d}^s\in\C{F} (\m{x} (C))\cap\C{D}$ we have
$\m{x} (C)\tr\m{Ad}^s \le 0$. Once this is shown, the proof will be complete,
since by
Lemma \ref{lemmaEdgeDesc} $\forall\;\m{d}\in\C{F} (\m{x} (C))$
$\exists\; k\in\mathbb{N}_+$, $\alpha_1,\alpha_2,\ldots,\alpha_k\ge 0$,
$\m{d}^1,\m{d}^2,\ldots,\m{d}^k\subseteq\C{F} (\m{x} (C))\cap\C{D}$ such
that
$\m{d} = \sum_{s = 1}^k \alpha^{s} \m{d}^s$, and therefore
\[
\m{x} (C)\tr\m{Ad} =
\sum_{s\in [k]} \alpha^s \m{x} (C)\tr\m{A}\m{d}^s \le 0\;.
\]
So, suppose that $\m{d}^s\in\C{F} (\m{x} (C))\cap\C{D}$. Then
$\m{d}^s = (\m{e}_i - \m{e}_j)$ for some $i\neq j$.
Moreover,
by definition of $\C{F} (\m{x} (C))$ we have $x_j (C) > 0$; that is,
$j\in C$. Hence, 
\begin{eqnarray}
\m{x} (C)\tr\m{A}\m{d}^s & = & \sum_{t = 1}^n \sum_{r = 1}^n x_t (C) a_{tr} d^s_r\nonumber\\
& = & \frac{1}{|C|}\sum_{t\in C} (a_{ti} - a_{tj})\nonumber\\
& = & \frac{1}{|C|}\bigg[\sum_{t\in C\backslash\{j\}} (a_{ti} - 1)
+ a_{ji}\bigg]\;.
\label{sinceCclique}
\end{eqnarray}
Here, 
(\ref{sinceCclique})
follows from the fact that $C$ is a clique, $j\in C$, and $a_{jj} = 0$.
We now consider two cases.

{\bf Case 1:} $i\in C$.
In this case, since $a_{ii} = 0$,
the right hand side of (\ref{sinceCclique}) equals
\begin{equation}\label{eq22}
\frac{1}{|C|}\bigg[\sum_{t\in C\backslash\{i,j\}} (a_{ti} - 1) +
a_{ji} - 1\bigg]\;.
\end{equation}
But since $i\in C$ and $C$ is a clique, we have that
$a_{ti} = 1$ for every $t\in C\backslash\{i,j\}$. Hence, (\ref{eq22}) is equal
to
\[
\frac{1}{|C|} (a_{ji} - 1)\le 0\;.
\]
Thus, $\m{x} (C)\tr\m{A}\m{d}^s \le 0$.

{\bf Case 2:} $i\notin C$.
In this case, we have that
(\ref{sinceCclique}) is equal to
\begin{equation}\label{eq23}
\frac{1}{|C|}\bigg[\sum_{t\in C\backslash\{i,j\}} (a_{ti} - 1)
+ a_{ji}\bigg]\;.
\end{equation}
But since $C$ is maximal and $i\notin C$ there must exist some $k\in C$
such that $a_{ik} = 0$. If $k = j$, then (\ref{eq23}) is less than or equal
to zero, since $a_{ji} = 0$ and
each of the terms in the first summation is less than or
equal to zero. On the other hand, if $k\neq j$, then there exists a term
in the first summation of (\ref{eq23}) which is equal to $-1$, and hence, since
$a_{ji} \le 1$, (\ref{eq23}) is less than or equal to zero.
Thus, $\m{x} (C)\tr\m{A}\m{d}^s \le 0$.
This completes the proof.
\endproof
\smallskip

\begin{lemma}\label{lemmaPhi}
Let
$\emptyset\neq S\subseteq\C{V}$. Then 
\begin{equation}\label{eq:nablaPhid}
\nabla \Phi (\m{x} (S)) \m{d} \le 0 \quad \forall\;\m{d}\in\C{F}(\m{x} (S))\;.
\end{equation}
Moreover,
if
$\nabla^2 \Phi (\tilde{\m{x}}) \succ \m{0}\;\forall\;\tilde{\m{x}}\in\Delta$,
then
the inequality in {\rm (\ref{eq:nablaPhid})}
is strict whenever
$\Supp{\m{d}}\not\subseteq S$.
\end{lemma}
\smallskip

\proof
First, observe that if (\ref{eq:nablaPhid}) holds in the case
when $\nabla^2 \Phi (\tilde{\m{x}})\succ\m{0}\;\forall\;\tilde{\m{x}}\in\Delta$,
then it also holds for any $\Phi$
(satisfying {\rm(C1)} -- {\rm(C3)}).
The argument is as follows: If $\Phi$ satisfies {\rm(C1)} -- {\rm(C3)}, then
for all sufficiently large $k\in\mathbb{N}$, the regularization
function $\Phi^{(k)} := \Phi + \frac{1}{k} \|\cdot\|_2^2$ also satisfies
{\rm(C1)} -- {\rm(C3)},
and moreover for any $\tilde{\m{x}}\in\Delta$ we have that
\[
\nabla^2 \Phi^{(k)} (\tilde{\m{x}}) = \nabla^2\Phi (\tilde{\m{x}}) +
\frac{2}{k}\m{I} \succ \m{0}\;.
\]
Hence,
\begin{equation}\label{Phikd}
\nabla \Phi (\m{x} (S)) \m{d} + \frac{2}{k}\m{x} (S)\tr\m{d}
=
\nabla \Phi^{(k)} (\m{x} (S)) \m{d} \le 0
\quad
\forall\;\m{d}\in\C{F} (\m{x} (S))\;.
\end{equation}
But then, taking the limit of {\rm (\ref{Phikd})}
as $k\rightarrow\infty$ we obtain
{\rm(\ref{eq:nablaPhid})}.

So suppose that
$\nabla^2 \Phi (\tilde{\m{x}}) \succ \m{0}$ for every
$\tilde{\m{x}}\in \Delta$.
We first prove the following result for any $\m{x}\in\Delta$:
\begin{equation}\label{strongerresult}
\nabla \Phi (\m{x})(\m{e}_i - \m{e}_j) < 0\quad\forall\; i,j\in [n]\mbox{ such that } x_i < x_j\;.
\end{equation}
To see this,
let $\m{x}\in\Delta$,
suppose that $i,j\in [n]$ are such that $x_i < x_j$, and
let $t := x_j - x_i > 0$ and
$\bar{\m{x}} := \m{x} + t (\m{e}_i - \m{e}_j)$. Then
$\bar{\m{x}}\in\C{P}(\m{x})$, since
$\bar{x}_i = x_j$, $\bar{x}_j = x_i$, and $\bar{x}_k = x_k\;\forall\;k\neq i,j$.
So, by (C3) we have
$\Phi (\bar{\m{x}}) = \Phi (\m{x})$. But taking a Taylor expansion of
$\Phi$ about $\m{x}$, we have that for some
$\g{\xi}\in [\m{x},\bar{\m{x}}]\subseteq\Delta$
\[
\Phi (\m{x}) = \Phi (\bar{\m{x}}) = \Phi (\m{x}) +
t \nabla \Phi (\m{x}) (\m{e}_i - \m{e}_j) +
\frac{t^2}{2}(\m{e}_i - \m{e}_j)\tr\nabla^2 \Phi(\g{\xi})(\m{e}_i - \m{e}_j)\;,
\]
and hence
\[
0 = t \nabla \Phi (\m{x}) (\m{e}_i - \m{e}_j) +
\frac{t^2}{2}(\m{e}_i - \m{e}_j)\tr\nabla^2\Phi (\g{\xi})(\m{e}_i - \m{e}_j)\;.
\]
So since $\nabla^2\Phi (\g{\xi})\succ \m{0}$ (by assumption)
and $t > 0$, we have that
\[
\nabla \Phi (\m{x}) (\m{e}_i - \m{e}_j) < 0\;.
\]
Thus, (\ref{strongerresult}) is proved. Moreover, note that
by continuity of $\nabla\Phi (\cdot) (\m{e}_i - \m{e}_j)$
we have in addition that
\begin{equation}\label{strongerresult2}
\nabla \Phi (\m{x})(\m{e}_i - \m{e}_j) \le 0\quad\forall\; i,j\in [n]\mbox{ such that } x_i = x_j > 0\;.
\end{equation}
In fact, by symmetry it is easy to see that the
inequality in (\ref{strongerresult2}) must actually be an equality.

Now we prove the lemma.
Let $\emptyset\neq S\subseteq [n]$ and
let $\m{d}\in\C{F} (\m{x} (S))$. Then by Lemma
\ref{lemmaEdgeDesc}
there exist
$k\in\mathbb{N}_{+}$,
$\m{d}^1,\ldots,\m{d}^k\in\C{F}(\m{x} (S))\cap\C{D}$,
and $\alpha^1,\alpha^2,\ldots,\alpha^k\ge 0$ such that
$\m{d} = \sum_{s = 1}^k \alpha^s \m{d}^s$. Hence,
\begin{equation}\label{Phid}
\nabla \Phi (\m{x} (S))\m{d} = \sum_{s = 1}^k \alpha^s \nabla \Phi (\m{x} (S))\m{d}^s\;.
\end{equation}
Next, note that by definition of $\C{D}$
for each $s \in [k]$ there exist $i,j\in [n]$ such
that $\m{d}^s = (\m{e}_i - \m{e}_j)$; moreover, since
$(\m{e}_i - \m{e}_j)\in\C{F}(\m{x} (S))$,
we have that $x_i \le \frac{1}{|S|} = x_j$. Hence, by
(\ref{strongerresult}), (\ref{strongerresult2}), and (\ref{Phid})
we have that $\nabla \Phi (\m{x} (S)) \m{d}^s \le 0$ for each $s \in [k]$.
So, by (\ref{Phid}) we have
$\nabla\Phi (\m{x} (S))\m{d} =
\sum_{s = 1}^k \alpha^s\nabla\Phi (\m{x} (S))\m{d}^s \le 0$.
Moreover,
in the case where $\Supp{\m{d}}\not\subseteq S$, there must exist some
$s\in [k]$ and some $i,j\in [n]$ such that $\m{d}^s = (\m{e}_i - \m{e}_j)$ with 
$i\notin S$ and $\alpha^s > 0$. By (\ref{strongerresult}), this implies
$\nabla\Phi (\m{x} (S))\m{d}^s < 0$. Hence, $\nabla\Phi (\m{x} (S))\m{d} < 0$.
This completes the proof.
\endproof
\smallskip

In the proof of the next proposition, we will use the following
well-known second-order sufficient optimality condition 
(see \cite{HagerHungerford13}):
A point $\m{x}\in\Delta$ is a local maximizer of (\ref{PMS}) if
\begin{equation}\label{firstordercond}
\nabla f (\m{x}) \m{d} \le 0\quad\forall\;\m{d}\in\C{F} (\m{x})\;,
\end{equation}
and
\begin{equation}\label{secondordercond}
\m{d}\tr\nabla^2 f (\m{x})\m{d} < 0\quad\forall\;
\m{0}\neq\m{d}\in\C{C} (\m{x})\;,
\end{equation}
where $\C{C} (\m{x})$ is the critical cone at $\m{x}$
defined by
\[
\C{C} (\m{x}) := \{\m{d}\in\C{F} (\m{x})\;:\; \nabla f (\m{x})\m{d} = 0\}\;.
\]

\begin{proposition}\label{propLocal}
Suppose that
$\nabla^2 \Phi (\tilde{\m{x}}) \succ \m{0}$ for every
$\tilde{\m{x}}\in \Delta$.
Then a point
$\m{x}\in\Delta$ is a
local maximizer of {\rm (\ref{PMS})}
if and only if
$\m{x} = \m{x} (C)$ for some maximal clique $C$. Moreover, every
local maximizer of {\rm (\ref{PMS})} is strict.
\end{proposition}
\smallskip

\proof
Suppose that
$\nabla^2 \Phi (\tilde{\m{x}}) \succ \m{0}$ for every
$\tilde{\m{x}}\in \Delta$.
First, we claim that every local maximizer of (\ref{PMS})
lies in $\Delta^0$.
To this end, let $\m{x}$ be any local maximizer of (\ref{PMS}).
If
$\Supp{\m{x}}$ is a clique, then clearly $\m{x}\in\Delta^0$. So,
suppose by way of contradiction that $\Supp{\m{x}}$ is not a clique.
By applying an argument similar
to the one given in the proof of Proposition \ref{propGlobal},
there exist indices
$i\neq j\in\Supp{\m{x}}$ such that for all
$t \in [-x_i, x_j]$ we have $\m{x} (t) = \m{x} + t (\m{e}_i - \m{e}_j)
\in\Delta$ and
$f (\m{x} (t)) > f (\m{x})$, where the strict inequality here
follows from (\ref{boilsdowntoPhi}) and the fact that
$\nabla^2 \Phi (\g{\xi}) \succ \m{0}$
for any $\g{\xi}\in [\m{x},\m{x} (t)]\subseteq\Delta$.
But this
contradicts the fact that $\m{x}$ is a local
maximizer of (\ref{PMS}).
Hence, $\Supp{\m{x}}$ is a clique and $\m{x}\in\Delta^0$.

So, every local maximizer of (\ref{PMS}) lies in $\Delta^0$ and
is therefore a local maximizer of (\ref{Delta0P}). Thus, by Proposition
\ref{lemmaDelta0P} every local maximizer of (\ref{PMS}) is equal to $\m{x} (C)$
for some maximal clique $C$. Since there are only finitely many
maximal cliques in $G$, there are only finitely many local maximizers of
(\ref{PMS}), which implies that
every local maximizer of (\ref{PMS}) is isolated, and is therefore
a strict local maximizer (see for instance \cite{Gibbons97}).


To complete the proof, we must show that for any
maximal clique $C$, $\m{x} (C)$ is a
local maximizer of (\ref{PMS}).
To see this, let $C$ be a maximal clique and let $\m{x} := \m{x} (C)$.
First, we show that the first-order condition (\ref{firstordercond}) holds.
Let $\m{d}^1\in\C{F} (\m{x})$. By Lemma \ref{firstorderLemma},
we have $\m{x}\tr\m{A}\m{d}^1
\le 0$, and by Lemma \ref{lemmaPhi},
$\nabla \Phi (\m{x}) \m{d}^1 \le 0$. Hence,
\begin{equation}\label{nablaf}
\nabla f (\m{x})\m{d} = 2 \m{x}\tr\m{Ad}^1 + \nabla\Phi (\m{x})\m{d}^1 \le 0\;.
\end{equation}
So, we will be done when we show that the second-order condition
(\ref{secondordercond}) holds. To see this, let $\m{0}\neq\m{d}^2\in \C{C} (\m{x})$
be arbitrary. Then,
\begin{equation}\label{nablaf=0}
0 = \nabla f (\m{x}) \m{d}^2 = 2\m{x}\tr\m{Ad}^2 + \nabla\Phi (\m{x})\m{d}^2
\le \nabla \Phi (\m{x})\m{d}^2\;,
\end{equation}
where the last inequality follows from Lemma \ref{firstorderLemma},
since $\m{d}^2\in\C{F} (\m{x})$.
Thus, by the second statement in Lemma \ref{lemmaPhi} we must have
$\Supp{\m{d}^2}\subseteq C$, implying that
$\m{x} + t \m{d}^2\in\Delta (C)\subseteq\Delta^0$ for all sufficiently small
$t > 0$.
And so by Part 2 of Lemma \ref{lemma1} we have
\[
{\m{d}^2}\tr\nabla^2 f (\m{x})\m{d}^2 < 0\;.
\]
This completes the proof.
\endproof

\section{Preliminary numerical results.}
\label{sectResults}

In this section, we conduct some preliminary numerical experiments on three
different regularization functions satisfying the
conditions outlined in Section \ref{sectGenReg} in order to give an
indication of the potential impact
of different regularization terms
on the performance of a local optimization
algorithm applied to (\ref{PMS}).
If in practice
a maximum clique (rather than merely a maximal clique) is sought,
the local optimization algorithm we employ in our experiments
would need to be incorporated
into a global optimization framework, such as branch and bound, in order to
ensure convergence to a global maximizer. 

\subsection{Regularization functions.}
\label{sectPenaltyFns}
We considered the following three regularization terms, with the
indicated choices of parameters:
\smallskip
\begin{eqnarray}
\Phi_{B} (\m{x}) &:=& \frac{1}{2}||\m{x}||_2^2,\label{eqnregs1}\\
\Phi_{1} (\m{x})\;&:=& \alpha_1\|\m{x}+\epsilon\m{1}\|_p^p ,\;\epsilon > 0,\ p > 2, \; \textstyle 0 < \alpha_1 < \frac{2}{p(p-1)(1+\epsilon)^{p - 2}} ,\label{eqnregs2}\\
\Phi_{2} (\m{x})\; &:=&\alpha_2 \sum_{i=1}^n (e^{-\beta x_i}-1),\; \beta > 0,\; 0 < \alpha_2 < \textstyle\frac{2}{\beta^2}.\label{eqnregs3}
\end{eqnarray}
\smallskip

%

\noindent
Here,
$\Phi_{B}$ is the 2-norm regularization function introduced by Bomze et al.
\cite{bomze1997evolution}, and
$\Phi_1$ is a generalization of $\Phi_B$ to
p-norms where $p > 2$.
$\Phi_2$ is a well-known (for instance, see \cite{BradleyMangasarianRosen})
approximation of the following non-smooth function:
\begin{equation}\label{zeronorm}
 \tilde \Phi(\m{x})=-\alpha_2 \|\m{x}\|_0\;,
\end{equation}
where $\|\m{x}\|_0 = \Supp{\m{x}}$.
The motivation behind the choice of $\Phi_2$ is as follows.
By definition of $\Phi_2$, maximizing
$\m{x}\tr\m{A}\m{x} + \Phi_2$ over $\Delta$ is closely
related to the problem of finding a solution to (\ref{MSQP}) which has
the smallest support (ie. the maximum sparsity). Following the argument
laid out in the proof of Proposition \ref{propGlobal}, from any global
maximizer of (\ref{MSQP}) which is \emph{not} a characteristic vector
for a maximum clique, there exists a path leading to another global maximizer
which \emph{is} a characteristic vector for a maximum clique and whose
support is strictly smaller than that of the starting point. Hence, the
global maximizers of (\ref{MSQP}) which have the smallest support are
necessarily the characteristic vectors for maximum cliques. Thus, 
$\Phi_2$ is a somewhat natural choice in our present context.

Next, we show that each of the regularization functions above satisfies
the conditions of Section \ref{sectGenReg}.
For any $\delta\in (0,\epsilon)$, where $\epsilon > 0$ is the
value used in the definition of $\Phi_1$, let
$X := conv (\cup_{i = 1}^n \C{B}_{\delta} (\m{e}_i)) \subseteq \mathbb{R}^n$,
where $\C{B}_{\delta} (\m{e}_i) =
\{\m{x}\in\mathbb{R}^n\;:\;\|\m{x} - \m{e}_i\|_2 < \delta\}$.
Then, $X$ is open, $X \supset \Delta$, and $\Phi_B$, $\Phi_1$, and
$\Phi_2$ are each well-defined on $X$. Moreover, it is easy to check that
$\Phi_B$, $\Phi_1$, and $\Phi_2$ are each twice continuously differentiable
over $X$ and that for any $\m{x}\in\Delta$ we have
\begin{eqnarray}
\nabla^2 \Phi_B (\m{x}) & = & \m{I}\succ \m{0}\;,\nonumber\\
\nabla^2 \Phi_1 (\m{x}) & = & \alpha_1 p (p - 1) Diag (\{(x_i + \epsilon)^{p - 2}\}_{i = 1}^n) \succ \m{0}\;,\; \mbox{and}\label{hessians}\\
\nabla^2 \Phi_2 (\m{x}) & = & \alpha_2 \beta^2 Diag (\{e^{-\beta x_i}\}_{i=1}^n)\succ \m{0}\;,\nonumber
\end{eqnarray}
where the positive definiteness of the Hessians in (\ref{hessians}) follows
from the choice of parameters. Hence, $\Phi_B$, $\Phi_1$, and $\Phi_2$ each
satisfy {\rm(C1)} strictly. Next, observe that
for any $\m{x}\in\Delta$ we have
\begin{eqnarray*}
\|\nabla^2 \Phi_B (\m{x})\|_2 &\; =\; & 1 < 2\;,\\
\|\nabla^2 \Phi_1 (\m{x})\|_2 &\; =\; & \alpha_1 p (p - 1) \max\;\{(x_i + \epsilon)^{p - 2}\}_{i = 1}^n\\
&\; <\; & \frac{2}{(1 + \epsilon)^{p - 2}} \max\;\{(x_i + \epsilon)^{p - 2}\}_{i = 1}^n\\
&\; \le \; & \frac{2}{(1 + \epsilon)^{p - 2}}  (1 + \epsilon)^{p - 2} = 2\;,\;\;
\mbox{and}\\
\|\nabla^2 \Phi_2 (\m{x}) \|_2 &\; =\; &
\alpha_2 \beta^2 \max\;\{e^{-\beta x_i}\}_{i = 1}^n\\
&\; <\; & 2 \max\;\{e^{-\beta x_i}\}_{i = 1}^n\\
&\; \le \; & 2\;.
\end{eqnarray*}
Thus, {\rm(C2)} is satisfied for each of $\Phi_B$, $\Phi_1$, and $\Phi_2$.
That {\rm(C3)} holds follows easily from the fact that
$\Phi_B$, $\Phi_1$, and $\Phi_2$ are each separable and the
coefficients associated with the terms $x_i$ are independent of $i$.

\subsection{The testing set.}

In the experiments, we considered different families of widely used maximum-clique instances belonging to the DIMACS benchmark \cite{johnson1996cliques}: 
\begin{itemize}
 \item \textbf{C family}: Random graphs Cx.y, where x is the number of nodes and y the edge probability;
 \item \textbf{DSJC family}: Random graphs DSCJx$\_$y. 
Here again, x is the number of nodes and y the edge probability;
 \item \textbf{brock family}: Random graphs with cliques hidden among nodes that have a relatively low degree;
 \item \textbf{gen family}: Artificially generated graphs with large, known embedded clique;
 \item \textbf{hamming family}: \mbox{hamminga-b} are graphs on a-bit words with an edge 
   if and only if the two words are at least hamming distance b apart;
  \item \textbf{keller family}: Instances based on Keller's conjecture
\cite{keller1930luckenlose} on tilings using hypercubes;
   \item \textbf{p$\_$hat family}: Random graphs generated with the p-hat generator, which is a generalization of the classical 
   uniform random graph generator. Graphs generated with p-hat have wider node degree spread and larger cliques than uniform graphs.
\end{itemize}
In Table \ref{tabDIMACS}, we report the names of the instances used (\textbf{instance}), the best known solutions (\textbf{best known}),
the number of nodes and edges in the instances (\textbf{nodes} and \textbf{edges}), and the median and interquartile range related to the graph 
degrees (\textbf{median} and \textbf{iqr} in column \textbf{graph degrees})
as well as the
median and iqr of the degrees of the nodes lying in the best known solution
(\textbf{median} and \textbf{iqr} in column \textbf{best degrees}).

\begin{table}
\caption{DIMACS instances used in the tests.\label{tabDIMACS}}
{\tiny
{\begin{tabular}{|l|r|r|r|r|r|r|r|}
\hline
\multicolumn{4}{|c|}{}		&\multicolumn{2}{|c|}{\textbf{graph degrees}}		&\multicolumn{2}{|c|}{\textbf{best degrees}}	\\
	        \hline
\textbf{instance}&\textbf{best known}	&\textbf{	nodes}	&\textbf{	edges}	&\textbf{	median}	&\textbf{	iqr}	&\textbf{	median}	&\textbf{	iqr}	\\
\hline
\hline
C125.9	&	34	&	125	&	6963	&	112	&	5	&	114.5	&	4.75	\\
C250.9	&	44	&	250	&	27984	&	224	&	6	&	227	&	5	\\
C500.9	&	57	&	500	&	112332	&	449	&	9	&	455	&	9	\\
C1000.9	&	68	&	1000	&	450079	&	900	&	13	&	907	&	11.25	\\
C2000.5	&	16	&	2000	&	999836	&	999	&	30	&	1006	&	11.5	\\
C2000.9	&	80	&	2000	&	1799532	&	1800	&	18	&	1803	&	15.25	\\
\hline
\hline
DSJC500$\_$5	&	13	&	500	&	125248	&	250	&	16	&	259	&	14	\\
DSJC1000$\_$5	&	15	&	1000	&	499652	&	500	&	20	&	503	&	23	\\
\hline
\hline
brock200$\_$2	&	12	&	200	&	9876	&	99	&	10	&	101	&	11	\\
brock200$\_$4	&	17	&	200	&	13089	&	131	&	8	&	134	&	6	\\
brock400$\_$2	&	29	&	400	&	59786	&	299	&	10	&	299	&	9	\\
brock400$\_$4	&	33	&	400	&	59765	&	299	&	11	&	299	&	9	\\
brock800$\_$2	&	24	&	800	&	208166	&	521	&	18	&	516.5	&	20.25	\\
brock800$\_$4	&	26	&	800	&	207643	&	519	&	18.25	&	512	&	20.25	\\
\hline
\hline
gen200$\_$p0.9$\_$44	&	44	&	200	&	17910	&	180	&	8	&	179.5	&	4.25	\\
gen200$\_$p0.9$\_$55	&	55	&	200	&	17910	&	179	&	7.25	&	179	&	5.5	\\
gen400$\_$p0.9$\_$55	&	55	&	400	&	71820	&	360	&	13.25	&	359	&	6	\\
gen400$\_$p0.9$\_$65	&	65	&	400	&	71820	&	361	&	14	&	359	&	9	\\
gen400$\_$p0.9$\_$75	&	75	&	400	&	71820	&	359	&	13	&	359	&	8	\\
\hline
\hline
hamming8-4	&	16	&	256	&	20864	&	163	&	0	&	163	&	0	\\
hamming10-4	&	40	&	1024	&	434176	&	848	&	0	&	848	&	0	\\
\hline
\hline
keller4	&	11	&	171	&	9435	&	110	&	8	&	112	&	17	\\
keller5	&	27	&	776	&	225990	&	578	&	38	&	578	&	33	\\
keller6	&	59	&	3361	&	4619898	&	2724	&	50	&	2724	&	50	\\
\hline
\hline
p$\_$hat300-1	&	8	&	300	&	10933	&	73	&	39	&	103	&	20	\\
p$\_$hat300-2	&	25	&	300	&	21928	&	146.5	&	73	&	213	&	18	\\
p$\_$hat300-3	&	36	&	300	&	33390	&	224	&	38	&	251	&	15.25	\\
p$\_$hat700-1	&	11	&	700	&	60999	&	174.5	&	87	&	250	&	22.5	\\
p$\_$hat700-2	&	44	&	700	&	121728	&	353	&	177.5	&	508	&	31.5	\\
p$\_$hat700-3	&	62	&	700	&	183010	&	526	&	89	&	602	&	14	\\
p$\_$hat1500-1	&	12	&	1500	&	284923	&	383	&	197	&	509	&	82	\\
p$\_$hat1500-2	&	65	&	1500	&	568960	&	763	&	387	&	1100	&	37	\\
p$\_$hat1500-3	&	94	&	1500	&	847244	&	1132.5	&	192	&	1297.5	&	25.75	\\
\hline
\end{tabular}}}{}
\end{table}

\subsection{Experiments.}
To conduct our tests,
we developed a multistart framework in {\tt MATLAB} that uses a hybrid algorithm as local optimizer. It combines the {\tt fmincon} 
solver with the Frank-Wolfe method \cite{frank1956algorithm}.
For each instance,
we ran 100 trials each with a different randomly generated point
in $\Delta$ as a starting guess. The same starting guesses were used for
all formulations. Since the iterates of the Frank-Wolfe method we used
are only guaranteed to converge to a point satisfying the first-order
conditions, we omitted the trials in which the final iterate was not
a true local optimizer from the
statistical computations in the tables below.
(Another potential
way of dealing with this issue, which we leave to a future work, would be to
take a step in an ascent direction whenever the final iterate is not a local
maximizer, and then rerun the algorithm using the new point as a starting
guess.)
Once a local maximizer $\m{x}^*$ is obtained, the associated clique is
constructed by taking
$C^* := \Supp{\m{x}^*}$.
In our experiments, we used the parameter values
$p =3$, $\epsilon=10^{-9}$, and $\beta = 5$. Furthermore,  $\alpha_1$ and $\alpha_2$ were suitably chosen in order to satisfy condition (\ref{eqnregs2}) and (\ref{eqnregs3}) respectively.
All the tests were performed on an  Intel Core i7-3610QM 2.3 GHz, 8GB RAM.

In Table \ref{tabDIMACSres}, we report 
the largest clique size
obtained (\textbf{max}), mean (\textbf{mean}), standard deviation (\textbf{std}) 
and the average CPU time (\textbf{CPU time}) over the 100 runs
for each instance and each of the three formulations.
When the largest clique
obtained is the same as the largest known clique size,
the result is reported in bold.
Observe that with the exception of the 
\textbf{keller} instances, 
the average clique size obtained
from using either of $\Phi_1$ or $\Phi_2$
was strictly larger than the
average obtained from using $\Phi_B$.
Overall, $\Phi_1$ performed slightly better than $\Phi_2$
yielding a strictly larger clique in 17 out of the 33 instances.
Next,
taking a look at the results related to the
\textbf{C} and \textbf{p\_hat} families, we notice that as the number of nodes
increases (and also the number of edges increases) finding a solution close to the best known gets harder and harder for all three formulations.
This is likely due to the simplicity of the global optimization 
approach that we use to solve the problem.
However, for the smaller instances in these groups $\Phi_1$ and $\Phi_2$
performed quite well. In particular, the formulation using
$\Phi_1$ found the largest known clique size in 6 of the instances.
The \textbf{DSJC}, \textbf{brock}, and \textbf{gen} families
all confirm the good behavior of the proposed formulations.
Indeed, in all cases the solutions found were closer (sometimes
significantly) to 
the best known clique than the ones found using $\Phi_B$.

\begin{table}
\caption{Results obtained on the DIMACS instances.\label{tabDIMACSres}}
{
\resizebox{12cm}{!} {
\begin{tabular}{|l|r|r|r|r|r|r|r|r|r|r|r|r|}
\hline
 &\multicolumn{4}{|c|}{$\Phi_{B}$}&\multicolumn{4}{|c|}{$\Phi_{1}$}&\multicolumn{4}{|c|}{$\Phi_{2}$}\\
\hline
\textbf{instances}	&\textbf{max}	&\textbf{mean}	&\textbf{std}	&\textbf{CPU time}	&\textbf{max}	&\textbf{mean}	&\textbf{std}	&\textbf{CPU time}&\textbf{max}	&\textbf{mean}	&\textbf{std}	&\textbf{CPU time} \\
\hline
\hline
C125.9		&	\textbf{34}		&		32.83		&		0.92		&		0.13		&	\textbf{34}		&	33.17	&	0.38	&	0.69	&	\textbf{34}		&		33.22		&		0.42		&		0.25		\\
C250.9		&		40		&		37.08		&		1.30		&		0.43		&	\textbf{44}		&	40.79	&	1.09	&	1.31	&	\textbf{44}		&		40.77		&		1.13		&		0.77		\\
C500.9		&		49		&		45.51		&		1.46		&		1.80		&	54	&	51.03	&	1.57	&	5.14	&		54		&		50.92		&		1.68		&		3.36		\\
C1000.9		&		60		&		54.40		&		2.26		&		10.01		&	63	&	58.90	&	1.91	&	23.57	&		63		&		59.16		&		1.94		&		15.76		\\
C2000.5		&		13		&		11.46		&		0.70		&		56.09		&	15	&	12.90	&	1.08	&	101.40	&		15		&		12.95		&		0.94		&		66.69		\\
C2000.9		&		69		&		61.75		&		1.86		&		56.67		&	67	&	63.35	&	1.61	&	106.22	&		68		&		63.89		&		1.96		&		76.71		\\
\hline
\hline
DSJC500$\_$5		&		10		&		9.43		&		0.61		&		2.20		&	11	&	10.45	&	0.55	&	4.59	&		12		&		10.19		&		0.44		&		2.36		\\
DSJC1000$\_$5		&		13		&		10.68		&		0.97		&		10.29		&	14	&	11.86	&	0.79	&	25.53	&		14		&		11.53		&		1.20		&		13.69		\\
\hline
\hline
brock200$\_$2		&		9		&		8.00		&		0.58		&		0.25		&	10	&	9.30	&	0.67	&	0.64	&		10		&		9.04		&		0.21		&		0.37		\\
brock200$\_$4		&		14		&		12.60		&		0.87		&		0.30		&	15	&	13.44	&	0.63	&	1.35	&		15		&		13.40		&		0.64		&		0.43		\\
brock400$\_$2		&		23		&		19.84		&		1.00		&		1.16		&	24	&	22.14	&	1.26	&	3.49	&		24		&		21.80		&		1.23		&		1.60		\\
brock400$\_$4		&		22		&		19.77		&		1.33		&		1.10		&	24	&	21.47	&	1.15	&	2.62	&		24		&		21.46		&		1.05		&		1.51		\\
brock800$\_$2		&		18		&		15.60		&		1.12		&		6.14		&	20	&	17.32	&	0.96	&	15.58	&		20		&		17.28		&		0.94		&		7.81		\\
brock800$\_$4		&		17		&		15.22		&		0.94		&		6.53		&	20	&	17.12	&	1.05	&	17.79	&		19		&		17.08		&		1.00		&		8.83		\\
\hline
\hline
gen200$\_$p0.9$\_$44		&		36		&		33.37		&		1.25		&		0.31		&	40	&	37.51	&	1.05	&	1.45	&		40		&		37.43		&		1.12		&		0.49		\\
gen200$\_$p0.9$\_$55		&		40		&		36.93		&		1.02		&		0.31		&	41	&	38.95	&	1.14	&	1.06	&		41		&		38.98		&		1.20		&		0.41		\\
gen400$\_$p0.9$\_$55		&		48		&		44.19		&		1.57		&		0.96		&	51	&	49.01	&	0.87	&	3.22	&		51		&		49.03		&		0.86		&		1.64		\\
gen400$\_$p0.9$\_$65		&		48		&		43.97		&		2.28		&		1.13		&	51	&	48.05	&	1.75	&	4.14	&		51		&		48.25		&		1.90		&		1.79		\\
gen400$\_$p0.9$\_$75		&		46		&		42.93		&		1.41		&		1.05		&	49	&	47.76	&	1.04	&	5.69	&		50		&		47.65		&		1.38		&		2.91		\\
\hline
\hline
hamming8-4		&	\textbf{16}		&		13.61		&		2.27		&		0.43		&	\textbf{16}		&	15.73	&	1.00	&	0.66	&	\textbf{16}		&		15.73		&		1.01		&		0.51		\\
hamming10-4		&		34		&		30.54		&		1.23		&		10.69		&	\textbf{40}		&	33.45	&	1.37	&	35.23	&	\textbf{40}		&		33.47		&		1.42		&		18.99		\\
\hline
\hline
keller4		&		8		&		7.17		&		0.38		&		0.19		&	7	&	7.00	&	0.00	&	0.73	&		9		&		7.02		&		0.20		&		0.18		\\
keller5		&		16		&		15.02		&		0.15		&		4.40		&	15	&	15.00	&	0.00	&	7.54	&		15		&		15.00		&		0.00		&		4.39		\\
keller6		&		34		&		32.86		&		1.46		&		388.05		&	34	&	33.83	&	0.61	&	290.79	&		34		&		32.80		&		1.48		&		222.96		\\
\hline
\hline
p$\_$hat300-1		&		7		&		7.00		&		0.00		&		0.54		&	\textbf{8}		&	8.00	&	0.00	&	0.88	&	\textbf{8}		&		8.00		&		0.00		&		0.52		\\
p$\_$hat300-2		&		24		&		24.00		&		0.00		&		0.55		&	\textbf{25}	&	24.01	&	0.10	&	1.89	&		24		&		24.00		&		0.00		&		0.44		\\
p$\_$hat300-3		&		33		&		31.15		&		0.76		&		0.62		&	\textbf{36}		&	33.39	&	0.70	&	1.86	&	\textbf{36}		&		33.20		&		0.64		&		0.90		\\
p$\_$hat700-1		&		9		&		7.24		&		0.62		&		3.97		&	9	&	8.15	&	0.43	&	13.31	&		9		&		7.94		&		0.55		&		6.15		\\
p$\_$hat700-2		&		43		&		41.53		&		0.75		&		4.49		&	\textbf{44}		&	43.61	&	0.70	&	10.07	&	\textbf{44}		&		43.50		&		0.72		&		4.74		\\
p$\_$hat700-3		&		60		&		58.71		&		0.87		&		4.77		&	61	&	58.99	&	0.69	&	10.55	&		61		&		59.21		&		0.80		&		7.73		\\
p$\_$hat1500-1		&		11		&		9.00		&		0.86		&		23.26		&	11	&	9.36	&	1.12	&	41.27	&		11		&		9.83		&		0.65		&		25.40		\\
p$\_$hat1500-2		&		62		&		59.09		&		1.10		&		29.06		&	62	&	60.79	&	0.87	&	44.76	&		64		&		62.30		&		1.15		&		30.42		\\
p$\_$hat1500-3		&		87		&		81.90		&		1.80		&		31.76		&	88	&	86.99	&	0.39	&	225.12	&		92		&		88.21		&		1.98		&		38.19		\\
\hline
\end{tabular}}}{}
\end{table}

\section{Conclusions.}
\label{conclusions}
We described a general regularized continuous formulation for the MCP and
developed 
conditions which
guarantee an
equivalence between the original problem and the 
continuous reformulation in both a global and a local sense. We
have also proved the results
in a step by step manner
which we hope
reveals some of the underlying structural properties
of the formulation.
We further proposed two specific regularizers that satisfy the general 
conditions given in the paper, and compared the two related continuous formulations 
with the one proposed in \cite{bomze1997evolution} on different families of widely used maximum-clique 
instances belonging to the DIMACS benchmark. The numerical results, albeit still preliminary, 
seem to confirm the effectiveness of the proposed regularizers;
that is, when a local optimization method is applied
to the new regularized formulations, cliques of high quality can often
be obtained in reasonable computational times.


%
%
%

\section*{Acknowledgments.}
This work has been partly developed while the first author was an ER of the ``MINO: Mixed-Integer Nonlinear Optimization'' program funded by the European Union.
The authors therefore gratefully acknowledge the financial support of the
European
Union's Horizon 2020 Marie Curie Network for Initial Training (ITN) programme under grant agreement No. 316647.

\bibliographystyle{amsplain}
\bibliography{references} 
\end{document}